\input mumacros.tex
\centerline {\bf {\uppercase  {Milnor and finite type invariants of
 plat-closures}} }
\bigskip
\medskip
{\baselineskip=10pt{\msmall
\centerline{{\small E}{\msmall FSTRATIA } {\small K}{\msmall ALFAGIANNI\ 
{\msmall AND}\ \  {\small X}{\msmall IAO}-{\small S}{\msmall ONG}
{\small L}{\msmall IN }}}}}
%\bigskip
%\centerline {}
\vskip.40in
\centerline{\msmall To Joan Birman on the occasion of her 70th birthday}
%\p {{\small A}{\Small BSTRACT.}} {{ \msmall {
%We show that for an $n$-component,
%$n$-bridge link and a
%positive integer $m$, the following is true:
%If the longitudes of
%$L$ lie
%in the $(m+2)$-th term
%of the lower central series of the link group
%then all the  finite type invariants
%of orders $\leq m$ for $L$
%are the same as these of the $n$-component unlink. }}}
\vskip.40in
{\baselineskip=11pt \footnote { }{\p {\msmall 
The authors are supported in part by NSF.}}}
\p {\largebold Introduction}
\medskip
\medskip
Finite type invariants of knots or links
can be defined combinatorially
using only link projections in 
$S^3$. 
In this setting it can be seen
that every Jones-type polynomial invariant
(quantum invariants)
is equivalent to a sequence of finite type invariants.
See [B2, BN] and references therein.
Although Vassiliev's original approach to finite type
knot invariants ([V]) rests on topological foundations
it doesn't offer any insight into how the invariants
are related to the topology of the complement 
of an individual knot or link.

Milnor's ${\bar \mu}$-invariants  ([M1,2]) are 
integer link
concordance invariants
defined in terms of algebraic data extracted 
from the link group.
It is known that the ${\bar \mu}$-invariants
of length $k$  detect exactly when
the longitudes of the link
lie in the $k$-th term of the lower central series
of the link group.
As was shown in [BN2] and [Li]
(and further illustrated in [H-M]), these invariants can
be thought of as
the link homotopy 
(or link concordance)
counterpart of the theory of finite type invariants.
As it follows from the results in the above manuscripts
if all finite type invariants of orders $\leq m$
vanish for a link
$L$, then all its Milnor invariants of
length $\leq m+1$ vanish.

In [K-L] we related the finite type invariants of a knot
to some geometric properties of its complement.
We showed that the invariants of orders $\leq m$ detect when
a knot bounds a {\it regular} Seifert surface $S$, whose complement
looks, modulo the first 
${m+2}$ terms of the lower central series
of its fundamental group, like the complement
of a null-isotopy.

In the present note we use the techniques of [K-L]
to study finite type invariants
of $n$-component links
that admit an $n$-bridge 
presentation (see for example [Ro]).
These links are precisely the ones
realized as plat-closures of pure braids
([B1]).
We show that for these links the vanishing of finite type
invariants of bounded orders is equivalent
to the vanishing of the $\bar \mu$-invariants
of certain length.
More precisely we have:
\medskip
\p {\bf Theorem 1.} {\sl Let $L$ be an $n$-component,
$n$-bridge link
and 
let $U^n$ denote the $n$-component
trivial link.
Set $G_L=\pi_1 (S^3\setminus L)$ and let $m\in \N$.
Suppose that the longitudes of $L$
lie in the $(m+2)$th term
of the lower central series of $G_L$;
or equivalently that all its $\bar \mu$ invariants of
length $\leq m+2$ vanish.
Then all finite type invariants of orders $\leq m$
of $L$ are the same as these of $U^n$.}
\medskip
It is not known whether there exists a non-trivial
link
all of whose finite type invariants
are the same as these of the trivial
link with the same number of components.
As it was shown by T. Cochran
([C2]) the vanishing 
of all $\bar \mu$-invariants
for an $n$-component
and $n$-bridge link is equivalent
to the property of being trivial.
Thus Theorem 1 doesn't produce
examples
of non-trivial links with trivial finite type
invariants. However,
there may still be examples
of pure braid plat-closures
that have certain (quantum) polynomial invariants
the same as these of the trivial link.

A special class of
$n$-component, $n$-bridge
links is the class of  {\it pure links}, 
i.e. ordinary closures of pure braids.
For these links Theorem 1, in fact a stronger version of it; that
the vanishing of all $\bar\mu$-invariants of length
$\leq m+1$ implies the vanishing of all finite type invariants of order 
$\leq m$, follows 
from the results of [BN2]. We 
don't know whether such a stronger statement holds in general, but we 
will have a little more to say about this in Remark 3.5.

To explain the idea of the proof of Theorem 1
let us recall that a Seifert surface $S$ 
is 
{\it regular} if it has a spine $\Sigma$ whose embedding in $S^3$,
induced by the embedding $S\subset S^3$, is
isotopic to the standard embedding of a bouquet of circles. Such a spine
will be called a {\it regular spine} of $S$. 
In [K-L] we studied a special class of projections of a
regular surface $S$,  which allowed us to
find a correspondence between certain
algebraic properties in $\pi_1(S^3 \setminus S)$
and the topology of $S^3 \setminus S$.
We may realize an
$n$-component, $n$-bridge link $L$
as ``a half'' of a regular spine of a regular Seifert surface
for which the techniques of [K-L] apply.
In particular, if all the longitudes of $L$
lie in the $(m+2)$-th term
of the lower central series of $G_L$,
we conclude that $L$ may be trivialized 
in $2^{m+1}-1$ ways by changing crossings.
This in turn is enough to show that all the invariants of orders $\leq m$
are the same as these of $U^n$.
%The other direction of the theorem follows 
%from the results 
%in [BN2] or [Li]. 

We organize this paper as follows:
In $\S 1$ we recall some basic facts about 
finite type invariants of links.
In $\S 2$ we study links that are
plat-closures of pure braids
and we recall from [K-L]
the results about projections of Seifert surfaces
that we need for this paper.
In $\S 3$ we prove Theorem 1 and give a few corollaries. 

\medskip

\noindent{\it Acknowledgment:} We thank James Conant
for a careful reading of [K-L]. His questions made us aware
of an oversight in its original version. 

\medskip
\medskip
\p {\largebold 1. Finite type invariants of links}
\medskip
\smallskip
For details on finite type invariants of knots and links
and for an extended bibliography the reader is
referred to [B2] or [BN1]. Below we
will  briefly recall the definition and a few basic facts that 
we will need in this note.

A {\it singular link} $L\subset S^3$ is an immersion
$\coprod S^1 \longrightarrow S^3$
whose only singularities are finitely many transverse double points.

Let ${\cal L}_n$ be the rational vector space generated by the set of ambient 
isotopy
classes of oriented, singular links  with  exactly $n$ double points.
A link invariant  $f$ can be extended to an invariant 
of singular links
by defining
\medskip
\centerline{\epsfysize=.4in\epsfbox[49 56 566 735]{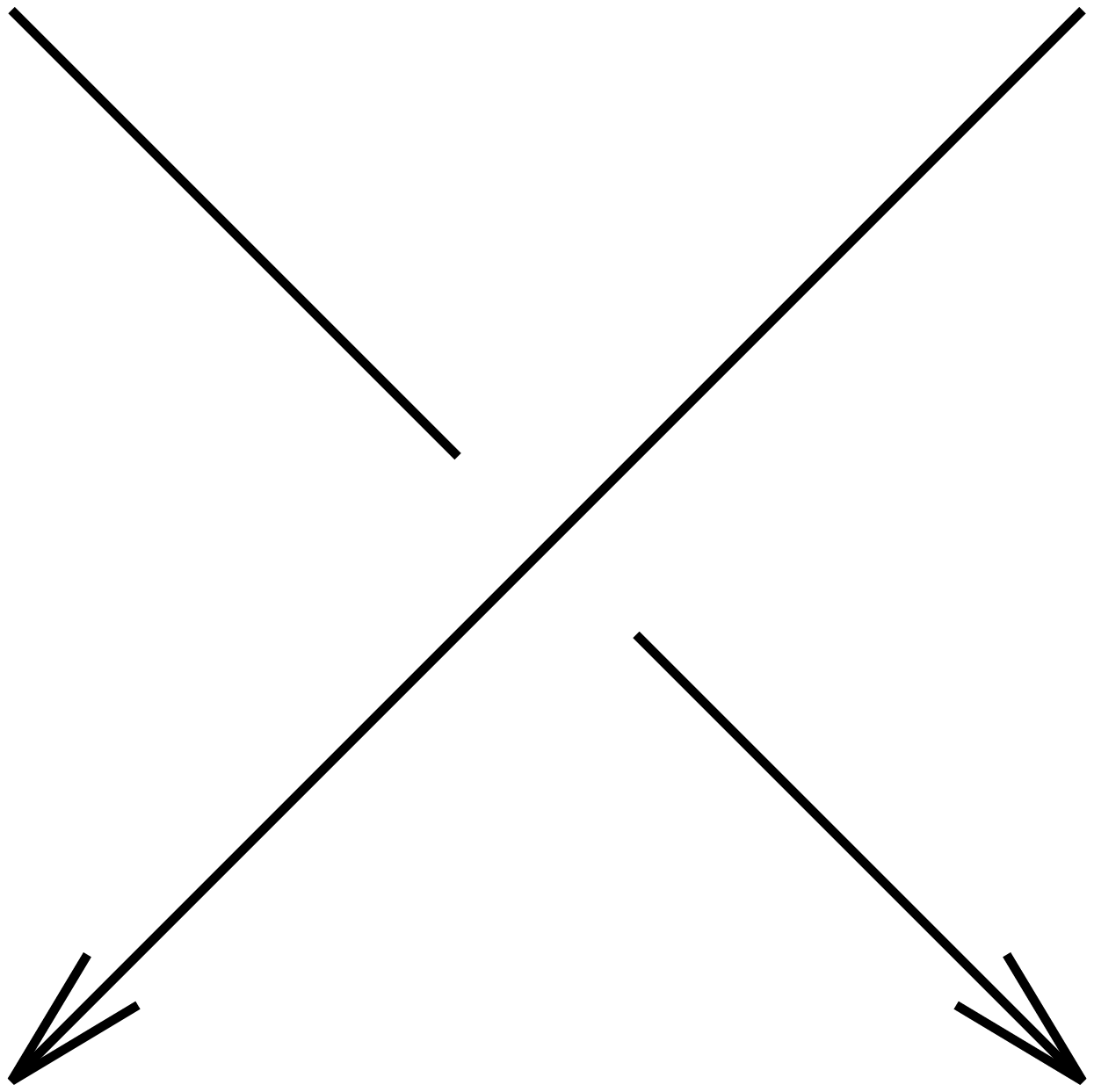}}
\medskip

\p for every triple of singular links which differ at one 
crossing as indicated.
Notice that  ${\cal L}_n$ can be 
viewed as a subspace of $\cal L$
for every $n$, by identifying any singular link in 
${\cal L}_n$
with the alternating sum of the $2^n$ links obtained by
resolving its double points.
Hence, we have a subspace filtration
$$\ldots \subset {\cal L}_n\ldots
\subset{\cal L}_2\subset {\cal L}_1 \subset {\cal L}$$

\medskip
\p {\bf Definition 1.1.} {\sl A  finite
type invariant of order $\leq n$
is a linear functional on the space 
${\cal L}/ {{\cal L}_{n+1}}$.
The invariants of orders  $\leq n$ form a subspace ${\cal F}_n$
of ${\cal L}^*$, the annihilator of the subspace
${\cal L}_{n+1}\subset {\cal L}$. We will say that
an invariant $f$ is of {\it order n} if  
$f$ lies in ${\cal F}_n$
but not in  
${\cal F}_{n-1}$. }
\medskip
To continue we need to introduce some notation.
Let $D=D(L)$ be a diagram of
a link  $L$, and let ${\cal C}={\cal C}(D)=
\{ C_1,\ldots, C_m\}$ be a collection of disjoint non-empty 
sets of crossings
of $D$. Let us denote by $2^{\cal C}$
the set of all subsets of $\cal C$. For an element $C\in 2^{\cal C}$ we will denote
by $D_C$ the link diagram obtained from $D$
by switching the crossings in all sets contained in $C$. 
So, all together, we
can get $2^m$ different link diagrams from 
the pair $(D,{\cal C})$. Notice that
each $C_i\in\cal C$ may contain more than one crossings.
\medskip
Following [G], two  $n$-component links $L_1$
and $L_2$ will be called 
{\it m-equivalent},
if $L_1$ has a link  diagram $D$ with the following property:
There exists ${\cal C}=\{C_1,\ldots, C_{m+1}\}$, a collection of $m+1$
disjoint
non-empty sets of crossings of $D$, such that $D_C$
is a diagram of $L_2$
for every non empty $C\in 2^{\cal C}$.
A link  $L$ which is {\it m-equivalent}
to the $n$-component unlink  will be called {\it m-trivial}.
\medskip
\p {\bf Proposition 1.2.} {\sl If $L_1$
and $L_2$ are  $m$-equivalent,
for some $m\in \N$,  then $f(L_1)=f(L_2)$
for all $f\in {\cal F}_m$.}
\medskip
{\it Proof.} Fix $m\in \N$ and let
$f\in {\cal F}_m$. 
Let ${\cal C}=\{C_1,\ldots, C_{m+1}\}$
be sets of crossings in a
projection of $L_1$,
satisfying the requirements
of $m$-equivalence.
By using
the skein relation of finite type invariants
first for all the crossings in $C_1$,
followed by these in $C_2$ and so on,
we obtain
$$f(L_1)= f(\Delta^{m+1})+
\sum_{i=1}^{m} f(\Delta^i)+ f(L_2),$$
\p where $\Delta^j$ is an element
in ${\cal L}_{j}$, for $j=1, \ldots, m+1$
and such that for $1\leq j\leq m$
$\Delta^j$ is identically zero in $\cal L$.
Then the conclusion follows 
since $f(\Delta^{m+1})=0$ by definition. \qed
\medskip
\medskip
\p {\largebold 2. Links in braid-plat form}
\medskip
\smallskip
In this section we will study links that are obtained as 
plat-closures of pure braids. 
\medskip
\smallskip
\p {\bf 2a. Preliminaries}
\medskip
Let $P_{2n}$ denote the pure braid group on $2n$ strings
(see for example [B1]). By a {\it pure $2n$-plat}
we will mean a pure braid $\sigma \in P_{2n}$
that is oriented in such a way so that any two successive
strings have opposite orientations.
We will say that an oriented link
$L=L(\sigma)$ is the closure of a pure $2n$-plat
if it is obtained by closing
the strings of a pure $2n$-plat as shown in Figure 1.
The following Lemma follows directly from the definitions.
\medskip
\p {\bf Lemma 2.1.} {\sl An oriented link $L$ is 
isotopic to the closure of a pure $2n$-plat
if and only if it is an $n$-component
$n$-bridge link.}
\medskip
\smallskip
\centerline{\epsfxsize=.9in\epsfbox[49 56 566 735]{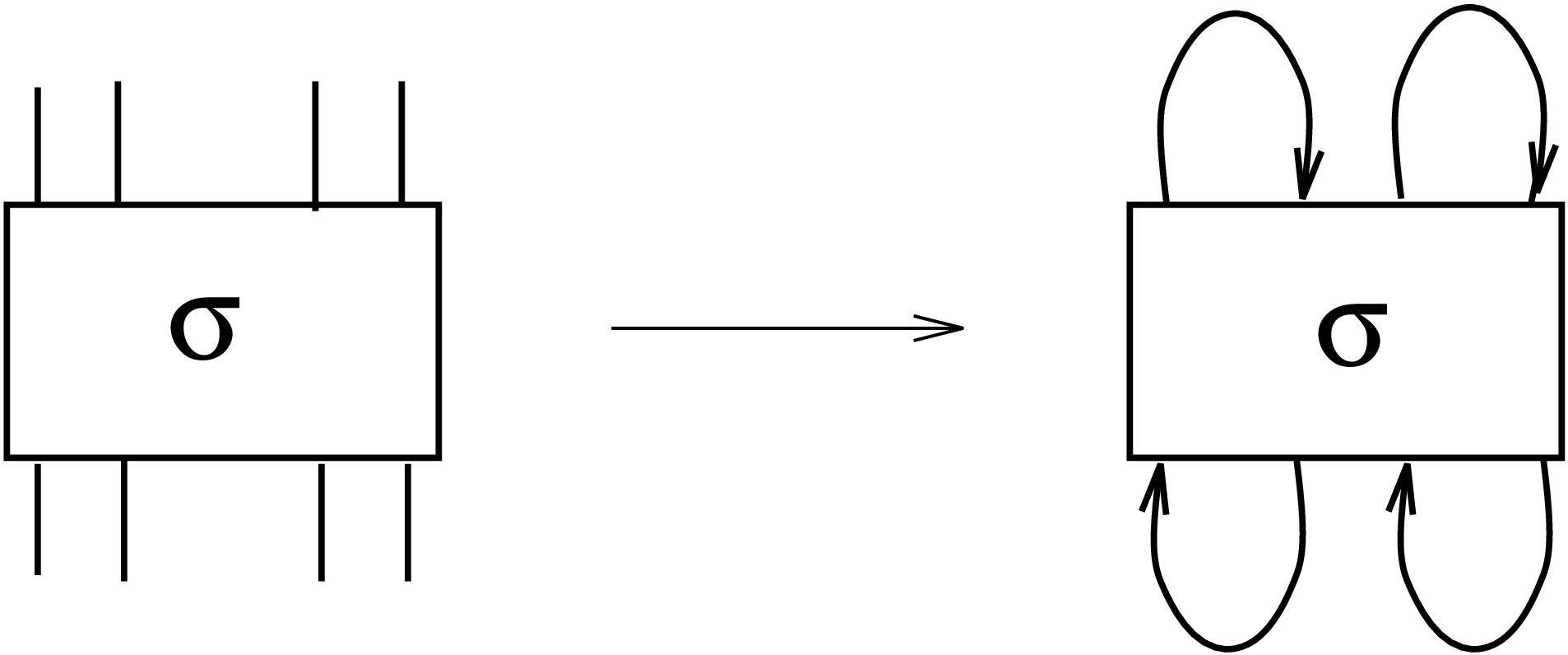} }

\centerline {{\bf Figure 1.} {\msmall A projection 
of a plat-closure of a pure braid.}}
\medskip
Let $L=L_1\coprod \ldots \coprod L_n$ be an $n$-component, 
$n$-bridge link and let $m_i$ and $l_i$ denote
the meridian and longitude of its $i$-th component,
respectively.
Fix a projection of $L$, on a projection plane $P\subset S^3$, 
as shown in Figure 1.
The meridians $m_1, \ldots, m_n$
can be realized by small linking circles 
$x_1, \ldots, x_n$, one for each component,  so that

i) lk$(x_i, L_j)=\delta_{ij}$; and

ii) their projections on the plane $P$ are simple curves 
disjoint from each other. 

Next we adjust the projection
so that each component realizes
the zero-framing.
This can be achieved by adding
a suitable number of kinks on the projection of each
component.
Now the longitude $l_i$,  is realized by a parallel
of $L_i$, for every $i=1, \ldots, n$.
 To continue, let $F=F(x_1, \ldots , x_n)$
denote  the free group on $x_1, \ldots , x_n$.
\medskip
\p {\bf Lemma 2.2.} ([C3, Theorem 3.1]) {\sl  Let $L$ be 
the closure of a pure $2n$-plat and let $x_1, \ldots , x_n$, $F$  and
$l_1, \ldots , l_n$
be as above. 
Set  $G_L=\pi_1(S^3 \setminus L)$ and
let $F \longrightarrow G_L$
be the map sending $x_i$ to the $i$-th meridian. 
Then $G_L$ has a presentation
$$<x_i\,;\,[x_i, W_i],\,i=1, \ldots, n>, \eqno (1)$$
\p where $W_i$ is a word 
in $F$ sent to $l_i$
under the meridian map above,
and $[x_i, W_i]=x_iW_ix_i^{-1}W_i^{-1}$.}
\medskip
For a group $G$ and a subgroup $H$
let $[G, H]$ denote the subgroup
of $G$ generated by commutators $ghg^{-1}h^{-1}$,
with $g\in G$ and $h \in H$. The lower central series
$\{ G^{(m)}\}_{m\in \N}$, of $G$ is defined
by
$G^{(1)}=G$ and
$$G^{(m+1)}=[G^{(m)}, G]$$
\p for $m\geq 1$.
\medskip
\p {\bf Lemma 2.3.}([C2, Theorem 4.1]) \ {\sl Let 
$L$ be an $n$-component, $n$-bridge
link and let
the notation be as in Lemma 2.2.
If  some longitude
$l_i$ lies in $G_L^{(m)}$ for some $m\in \N$, 
then
$W_i$ lies in $F^{(m)}$. In particular if
all the longitudes of $L$ lie in
$\cap_{m\in \N} G_L^{(m)}$, then $L$ is the trivial link.}
\medskip
{\it Proof.} For $m=1$ the conclusion is clearly true.
If $l_i$ lies in $G_L^{(m)}$ then $l_i$ lies in $G_L^{(m-1)}$
and thus $W_i \in F^{(m-1)}$, by induction.
By (1), $l_i$ must lie in the normal subgroup of 
$G_L$ generated by $F^{(m)}$ and the relations
$[x_i, W_i]$. Since $W_i \in F^{(m-1)}$
we obtain that $[x_i, W_i] \in F^{(m)} $.
Thus we conclude that
$W_i \in F^{(m)}$. 
Now the last part of the claim follows from the well known fact that
$\cap F^{(m)}=\{ 1\}$. \qed
\medskip
\smallskip
\p {\bf 2b. Plat-closures and spines of surfaces}
\medskip
 We will say that
a Seifert surface $S$ of a knot $K$ is 
{\it regular} if it has a spine $\Sigma$ whose embedding in $S^3$,
induced by an embedding $S\subset S^3$, is
isotopic to the standard embedding of a bouquet of circles. Such a spine
will be called a {\it regular spine} of $S$. 
In particular, $\pi=\pi_1(S^3 \setminus S)$
is a free group.
Let us start with an embedding $\Sigma \subset S\subset S^3$
as above.
For each circle in $\Sigma$, we may push it off 
the surface $S$ slightly in the normal directions. 
If these push-offs are all disjoint (this can be achieved by either pushing
circles off to different sides of $S$ or in different distances away from $S$),
we get a link. By Lemma 2.2 of [K-L]
we may find an embedding $S\subset S^3$ so that
each of these links
is a $2g$-component,  $2g$-bridge link,
where $g$ is the genus of $S$.

Now let $L$ be an $n$-component,
$n$-bridge link. We fix a projection, $p: L \longrightarrow P$, 
of $L$  as a closure of
a pure $2n$-plat, and such that the projection of each
component realizes the zero framing.
For every sequence
$${\bar \epsilon}:=(\epsilon_1, \ldots, \epsilon_{n}),$$
\p where each $\epsilon_i$ is equal  to $+1$ or $-1$,
we obtain a projection of a regular Seifert surface
$S^{\bar \epsilon}_P$, as follows:
Pick a point $Q$
on the projection plane and connect
it to a point on the projection of each $L_i$, by
an arc $\lambda_i$. This way we obtain a projection of
a bouquet of $n$-circles, say $W_n$.
For every circle in $W_n$ we add an unlinked loop $L_{i^{'}}$,
which contains a positive or a negative  kink
according to whether $\epsilon_i$ is equal  to $+1$ or $-1$ . 

This is done in such
a way so that the four arcs of $L_i$ and 
$L_{i^{'}}$ in a
disc neighborhood of $Q$,
appear in alternating order. See Figure 2. 
Call the resulting bouquet $W_{2n}$.
Now let $S^{\bar \epsilon}_p$
be the regular Seifert surface
{\it associated} to the given projection
of $W_{2n}$ [see $\S 2$ of K-L].
The boundary of $S^{\bar \epsilon}_p$
is clearly a knot. We shall assume further 
that the push-off of each $L_i$ 
that we are going to consider is 
geometrically unlinked with all the $L'_i$'s.
\medskip
\vskip .20in
\smallskip
\medskip
\input epsf
\centerline{\epsfysize=.4in\epsfbox[49 56 566 735]{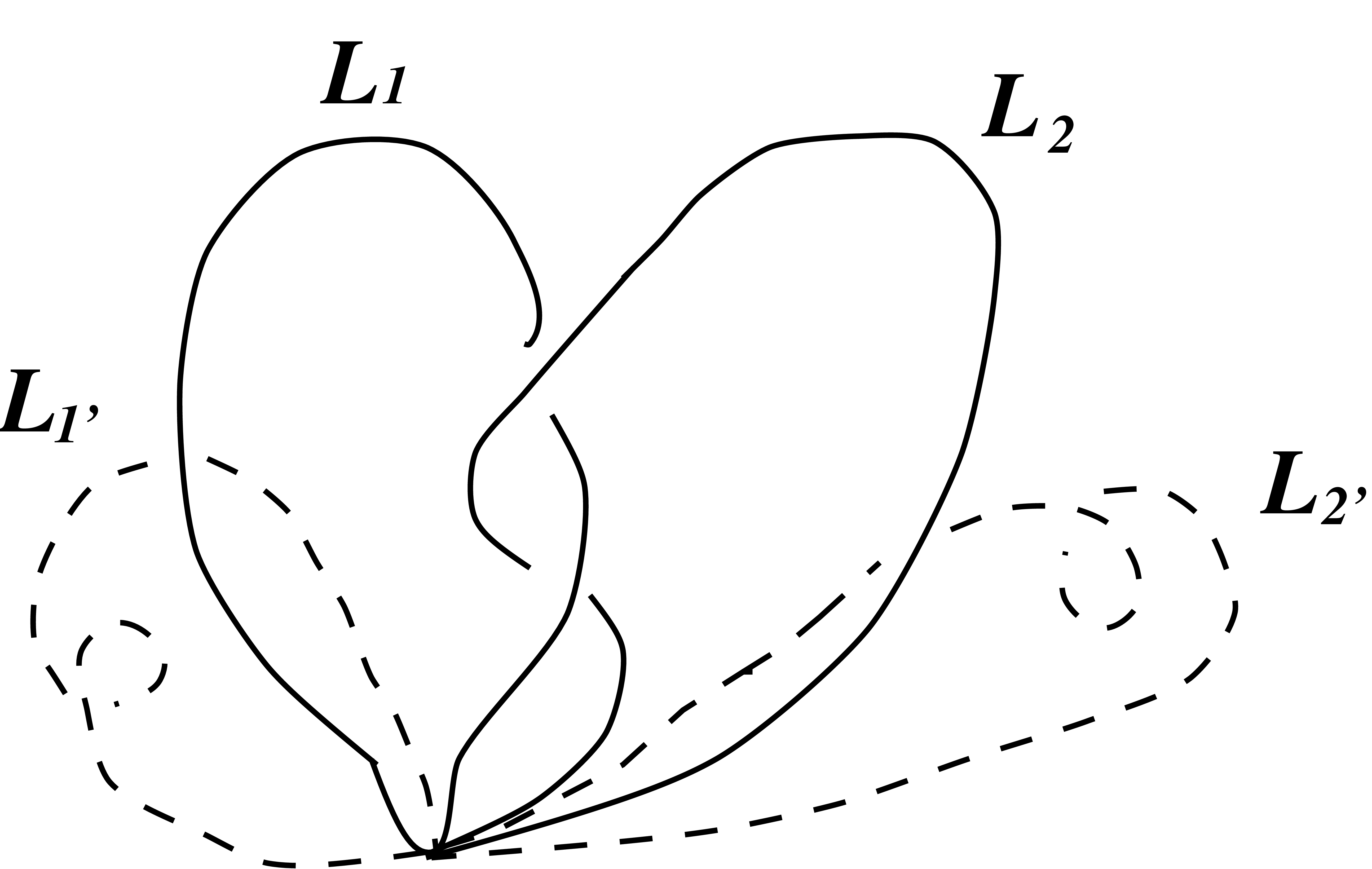} }
\bigskip
\medskip
\medskip
\medskip
\centerline{\bf Figure 2 }
\medskip
\smallskip
{\it Remark.} The construction of the surface 
$S_p^{\bar\epsilon}$ is not really essential
for the proof of our results here.
It only makes it easy 
for us to refer to results in [K-L].   
\medskip
\smallskip
\p {\bf 2c. Modifications of link projections}
\medskip
In this paragraph we will use the properties of projections of
Seifert surfaces established in [K-L]
to study projections of links that are closures of pure plats.
We begin with the following definition:
\medskip
\p{\bf Definition 2.4.} {\sl A {\it simple commutator} in
$F=F(x_1, \ldots , x_n)$ is a word in the form
of $[A,x^{\pm1}]$ or $[x^{\pm1},A]$ where $x$ is a generator and $A$ is a
simple commutator of shorter length. A {\it simple quasi-commutator} is a word
obtained from a simple commutator by finitely 
many insertions of
cancelling
pairs, $x_i^{\pm1}x_i^{\mp1}$ .}
\medskip
Any word in $F$ representing an element 
in $F^{(m)}$ can be changed
to a product of simple quasi-commutators of length 
$\geq m$ by finitely many
insertions of cancelling pairs.

Let
$L=L_1\cup \ldots \cup L_n$ be an $n$-component,
$n$-bridge link and set $G_L=\pi_1(S^3\setminus L)$.
Fix a projection, $p: L \longrightarrow P$,
on some plane $P$, and such that each component
realizes the zero-framing.
\medskip
\p {\bf Lemma 2.5.} \
{\sl Let $L$ and $p$ be as above and suppose that the longitude
$l=l_i$ corresponding to some  component $L_i$, lies in
$G_L^{(m)}$. Let $W=W_i$ be the word
of Lemma 2.3 for $l_i$, and let $W^{'}=c_1\ldots c_r$
be a word expressing $W$ as 
a product of simple quasi-commutators of
length $\geq m$.
There exists   
a projection $p_1:L\rightarrow P$ with 
the following properties:

\p (i) $p_1(L)$ is obtained from $p(L)$
by a finite sequence of Reidermeister 
moves 
on the component $L_i$; 

\p (ii) the word one reads out from $l$ 
(with respect to the new projection 
$p_1$), 
by 
picking up one letter for each
undercrossing of $L_i$ is $W^{'}$.}
\medskip
{\it Proof.} Fix a 
sequence $${\bar \epsilon}:=(\pm 1, \ldots, \pm 1)$$ and let 
$S=S^{\bar \epsilon}_p$,
be the regular Seifert surface
associated to the projection $p$ and to
${\bar \epsilon}$ as described in 2b.
Let $B_1, \ldots, B_n$ denote
the band of $S^{\bar \epsilon}_p$
whose core is realized by $L_1, \ldots, L_n$,
respectively. The set of  meridians
$x_1, \ldots, x_n$  described in 2a 
can be taken as part of a set of free generators of
$\pi_1 (S^3 \setminus S)$.
The projection $p$ gives rise to a projection $p(S)$
of $S$, so that one of the push-offs of $L_i$
is represented in $\pi_1 (S^3 \setminus S)$
by $W_i$.
Now Lemma 2.4 of [K-L] applies and we may
modify $p(S)$ to a projection ${\bar p}(S)$, 
so that the band $B_i$ corresponding to our fixed component $L_i$,
is in {\it good position with respect to $\bar p$}. See Definition
2.3 of [K-L]. Now we can apply Lemma 2.12
of [K-L] to the band $B_i$ and the projection
${\bar p}(S)$, to obtain a new
projection $p_1(S)$. By the proof of Lemma 2.12
of [K-L]  we see that $p_1$
if restricted to $L$ has all the desired properties. \qed
\medskip
We can decompose the component $L_i$ of Lemma 2.5 into 
sub-arcs
 $\delta_i$ and $\delta_i^{*}$  with
disjoint interiors and
such that: $\delta_i$
represents the word $W^{'}$
of Lemma 2.5 and $\delta_i^{*}$
is a small straight segment,
containing no crossings at all.

For a set of crossings $C$ in the projection
$p_1(L)$ we will denote
by $L_C$ 
the link obtained from $L$ by changing all 
the crossings
in $C$ simultaneously. Also, for $i=1, \ldots, n$,
let $L_i^C$ (resp. $\delta_i^C$) be the component 
(resp. the subarc of the component )
of $L_C$
corresponding to the component
(resp. the subarc  $\delta_i \subset L_i$ of the component)
$L_i$  of $L$.
\medskip
\p {\bf  Proposition 2.6.}
\ {\sl  Let $L$ 
be an $n$-component, 
$n$-bridge link and suppose that the longitude
$l=l_i$ corresponding to the component $L_i$ lies in
$G_L^{(m+1)}$. Also let
$p_1(L)$ be a projection  as in Lemma 2.5 for $L$, and let 
 $\delta_i$
and  $\delta_i^{*}$
be as above.
 We can find
a collection,
$ {\cal C}=\{ C_1, \ldots, C_m\} $, 
of sets of crossings on $p_1(L)$ such that

\p (i)  The crossings in each $C_i$ are crossings on $L_i$;

\p (ii) For every non-empty $C\in 2^{\cal {C}}$,
$L_C$ is an $n$-component, $n$-bridge link
and 

$\delta_i^C$ is isotopic to  $\delta_i^{*}$
in the
complement $S^3 \setminus L_C$.

\p (iii) The links $L_C$ are isotopic for all
 non-empty $C\in 2^{\cal {C}}$.

Here $2^{\cal {C}}$ is 
the set of all subsets of $\cal C$.}
\medskip
{\it Proof.} It follows 
from Proposition 3.3 of [K-L] and the discussion preceding it. \qed
\smallskip
\medskip
\p {\bf 2d. An Example}
\medskip
In Figure 3, we show a projection of a 2-component
Whitehead link $L$, as the closure of a pure
4-plat.
We have denoted the components of
$L$ by $A$ and $B$. Let $l_A$ and $l_B$ denote
the longitude of $A$ and $B$, respectively.
The fundamental group
$G_L=\pi_1(S^3\setminus L)$ is generated by 
the meridians $x$ and $y$. We have
$l_A \in G_L^{(3)}$.
In fact, from the 
Wirtinger presentation obtained from the given projection
we have $l_A=[[x, y], y^{-1}]$. 
\bigskip
\centerline{\epsfxsize=1in\epsfbox[49 56 566 735]{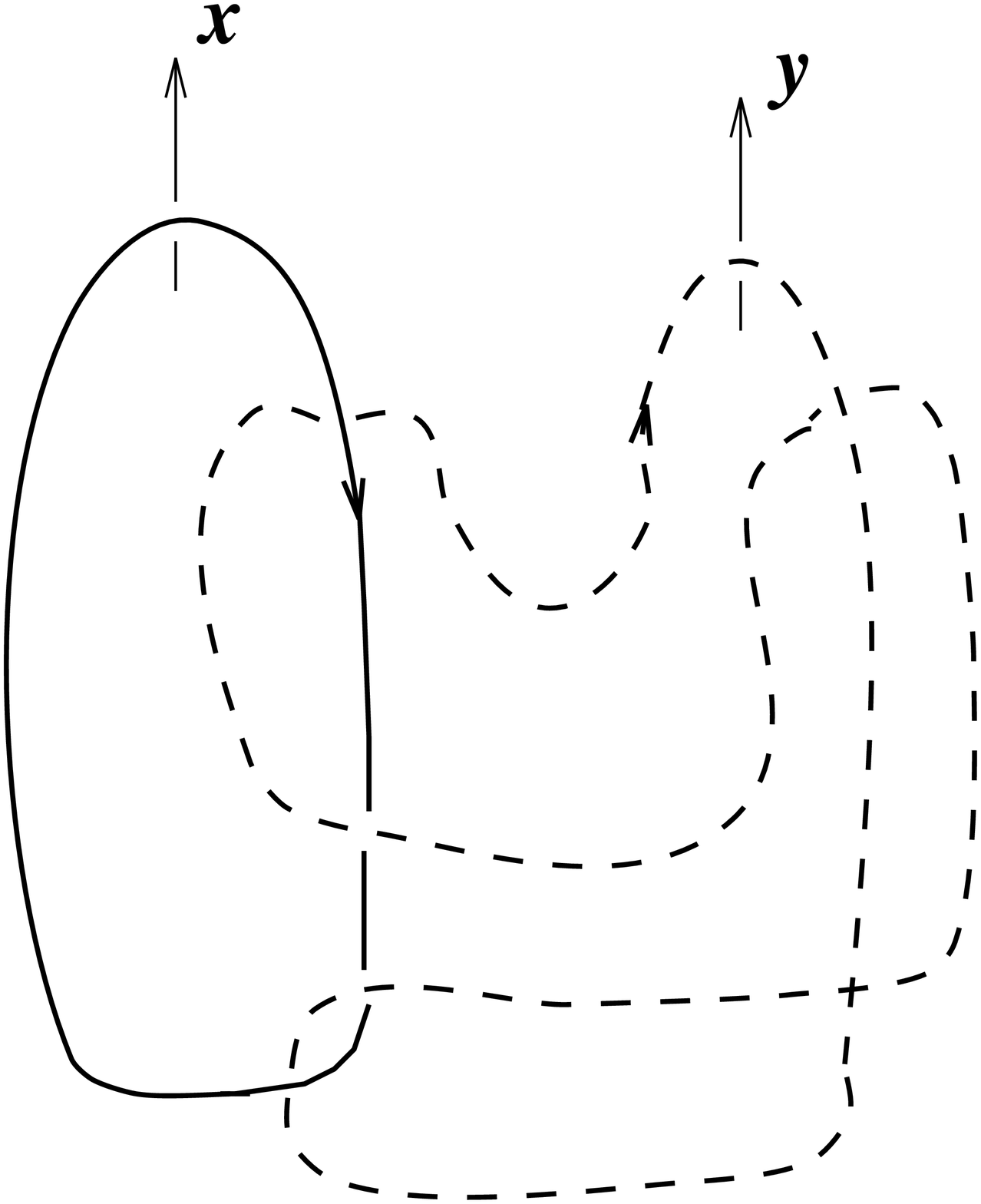} }
\medskip
\medskip
\centerline {{ \bf Figure 3.}\ {\msmall A projection of the Whitehead link. }}
\bigskip
In Figure 4 we have modified the projection
of $L$ so that it satisfies the properties
of Lemma 2.5. 
The solid (dashed, resp.)
line corresponds to the 
component $A$ ( resp. $B$)
of Figure 3. The word we read out when traveling along the solid line,
one letter for each undercrossing, is exactly
$W=[[x, y],y^{-1}]$. Notice that 
all the undercrossings of $A$ occur above the horizontal line $l$
shown in Figure 4.

\bigskip
\centerline{\epsfxsize=2.0in\epsfbox{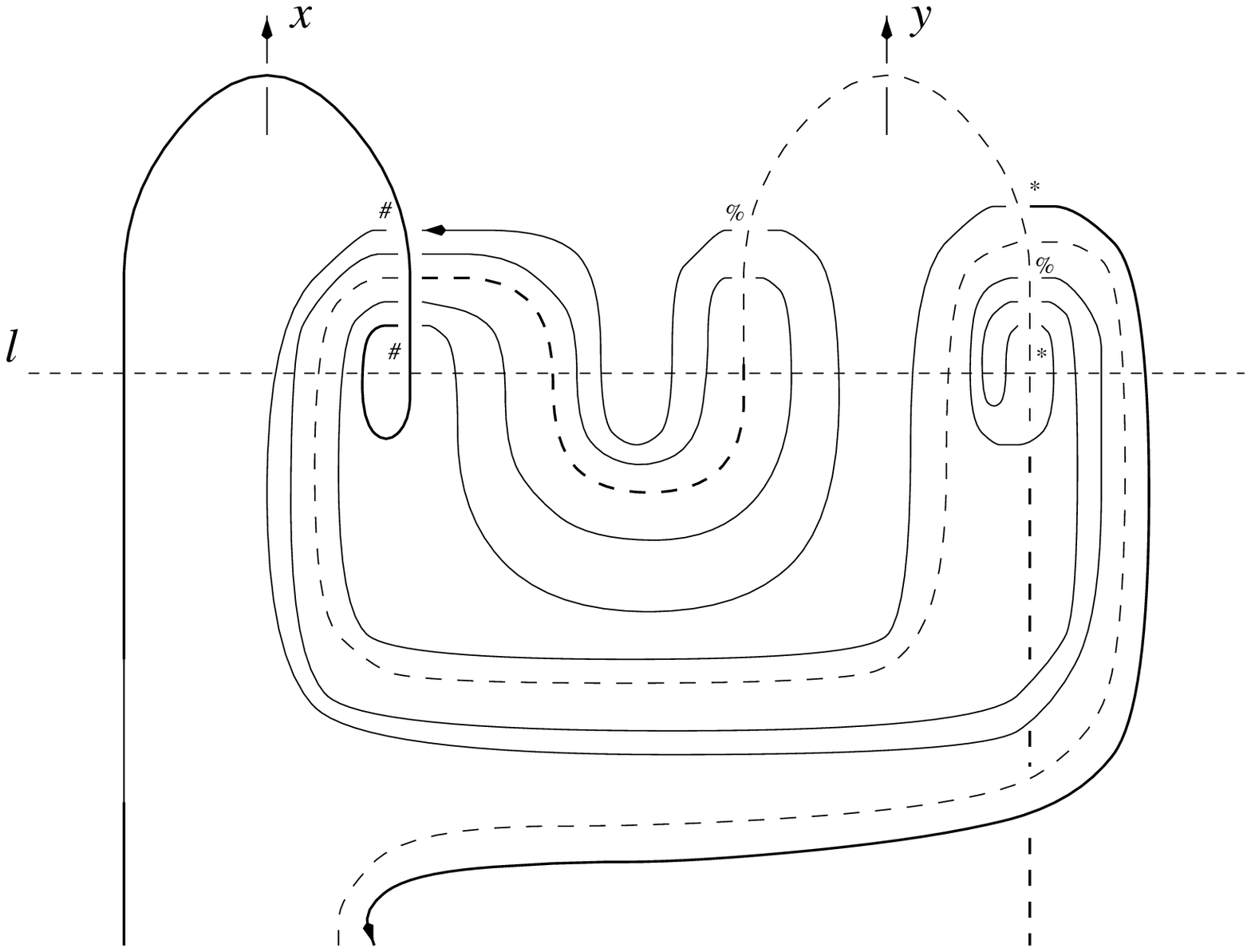} }
\medskip
\bigskip
\medskip
\centerline {{ \bf Figure 4.}\ {\msmall Realizing $l_A$
as a simple commutator of length three. }}
\medskip
The reader can verify that the two crossings marked by $``\%"$
on the left side of the $y$-hook and the two crossings marked
by $``\ast"$ can be used to trivialize the component
$A$ in three ways, and thus show that the link is 
{\it 1-trivial}. Of course, such sets of crossings can 
be spotted on the link diagram in Figure 3 as well. But our emphasis here
is the correspondence between such sets of crossings and the letters in the
word $[[x, y],y^{-1}]$. 
This example is meant to illustrate the correspondence between
the sets of crossings $C_i$ of Proposition 2.6 and the letters
in the words $W^{'}$ of Lemma 2.5, as well as
the proof of Theorem 3.1 below.

Let us remark that $L$ has a non trivial finite type invariant
of order three, coming from the Alexander-Conway polynomial.
There are no non-trivial $\bar \mu$-invariants of length five for 
2-component links.
As predicted by Theorem 3.1, $L$ also has a non-zero
$\bar \mu$-invariant,  of length four,  coming from the Sato-Levine invariant
(see for example [C1]). 

\vfill
\eject
\medskip
\medskip
\smallskip
\p {\largebold 3. Relations among link invariants }
\medskip
\smallskip
Our main purpose in  this section, is to show
that for $n$-component, $n$-bridge links
Milnor's $\bar \mu$-invariants control all
the finite type invariants.
 
Let us begin by briefly recalling
the basic idea and notation behind
Milnor's $\bar \mu$-invariants.
Let 
$L=\coprod_{i=1}^n L_i$ be an $n$-component link and 
suppose that the longitudes of $L$ lie in
$G_L^{(m)}$, 
for a positive integer $m$.
For a sequence $I=i_1, \ldots, i_{m+1}$
of numbers from the set $\{ 1, 2, \ldots, n\}$,
Milnor's ${\bar \mu}(I)$ invariant
is an integer concordance invariant.
The integer $m+1$ is called the length of ${\bar \mu}(I)$. 
If the longitudes of $L$ do not lie in
$G_L^{(m)}$ then the $\bar \mu$-
invariants of length $m+1$ are well defined only modulo
the ideal generated by certain invariants
of less length. For details the reader is referred to [M1, M2, C1].

The $\bar \mu$-invariants of length $m$ detect precisely when
the longitudes of $L$ lie 
in $G_L^{(m)}$. That is,  we have that Milnor's invariants 
of length $\leq m$ vanish for a link $L$ 
if and only if the longitudes of $L$ lie in $G_L^{(m)}$.
\medskip
\p {\bf Theorem 3.1.} {\sl Let $L$ be an $n$-component, $n$-bridge link
and let $U^n$
denote the $n$-component unlink.

a) If $f(L)=f(U^n)$
for all $f\in {\cal F}_m$, for some $m\in \N$, the
$\bar \mu$-invariants of $L$ of length  $\leq (m+1)$
vanish.
\smallskip 
b) If the $\bar \mu$-invariants of $L$ of length  $\leq (m+2)$
vanish for some $m\in \N$, then
$f(L)=f(U^n)$
for all $f\in {\cal F}_m$.}
\medskip
{\it Proof.} a) This direction follows  from 
the results in [BN2, Li]. For completness, we include
the proof here. 

Assume that $L$ is 
an $n$-component, link whose finite type
invariants of orders $\leq m$ ($m\in \N$)
 are the same as these
of the $n$-component unlink.
For $m=1$ the conclusion follows easily from the
fact that the linking numbers of $L$ are finite type 
of order $1$, and that the only Milnor invariants
of length $\leq 2$ come from linking numbers.
Inductively, suppose that the conclusion is
true for all integers $<m$, and let $L$ be an $n$-component
 link
whose finite type
invariants of orders $\leq m$ 
 are the same as these
of the $n$-component unlink.
By induction, all the $\bar \mu$-invariants of
length $\leq m$ vanish for $L$.
Thus, the longitudes of $L$
lie in $G_L^{(m)}$ and  the $\bar \mu$-invariants
of length  $m+1$ are well defined.
By the arguments in [BN2, Li]
they 
are of finite type of order $\leq m$.
Now the conclusion follows since the $\bar \mu$-invariants of 
an unlink vanish.

For the converse let $L$ be an $n$-component, $n$-bridge 
link with vanishing $\bar \mu$-invariants
of length $\leq (m+2)$. Then the longitudes
of $L$ lie in $G_L^{(m+2)}$. 
We will induct on the number of components of $L$
that are non-trivial in $S^3\setminus L$.
Suppose that  $L_i$ is such a component
and let $p_1(L)$ be  a projection of $L$,
that
has  the properties of Lemma 2.5.
By Proposition 2.6, $L$ is $m$-equivalent
to an an $n$-component and  $n$-bridge link $L_C$, such that

 (i) the component $L_i^C$ is null-isotopic
in the
complement $S^3 \setminus L_C$;

 (ii) the longitudes of $L_C$ lie in $\pi_1^{(m+2)}(S^3 \setminus L_C)$;
thus the $\bar \mu$-invariants
of length $\leq m+2$ vanish for $L_C$.

Now the conclusion follows by induction and Proposition 1.2. \qed
\medskip
By  Lemma 2.3 we see that for $n$-component, $n$-bridge
links the property of being trivial is equivalent to the vanishing
of all $\bar \mu $-invariants. Thus we have:
\medskip
\p {\bf Corollary 3.2.} {\sl Let $L$ be an
$n$-component, $n$-bridge link. We have
that $f(L)=f(U^n)$
for all $f\in {\cal F}_m$, and all $m\in \N$
if and only if $L$ is the $n$-component unlink.}
\medskip
A link $L$ is called {\it Brunnian} iff 
every proper sublink of $L$
is a trivial link.
\medskip
\p {\bf Corollary 3.3.}\ {\sl Let $L$
be an $n$-component, $n$-bridge
Brunnian link, with $n>3$. Then 
$f(L)=f(U^n)$ for all
$f\in {\cal F}_m$, and all $m < n-2$.}
\medskip
{\it Proof.} It follows from
Theorem 3.1 and the fact that Milnor's
invariants of length $\leq n-1$ vanish
for a Brunnian $n$-component link. \qed
\medskip
\p {\bf Remark 3.4.} 
It is an open problem whether there exist 
non-trivial links all
of whose finite type invariants are the same 
as these of an unlink.
In view of Corollary 3.2,
Theorem 3.1 doesn't produce
examples of such links.
On the other hand,
it is known that there exist links 
all of whose Milnor's invariants vanish
but they are not even concordant to the trivial
link ([C-O]). So the question is whether Theorem 3.1
can be extended to include any of these links
or whether the techniques of [K-L] can be extended
to surfaces whose complement is not a handlebody.
\medskip
\p {\bf Remark 3.5.} a)  Let $P_n$
be the pure braid group on $n$ strings
and $F=F(n)$ be the free group of rank $n$.
For a positive integer $m$, let
$${\cal A}_m: P_n \longrightarrow 
Aut({{F} / {F^{(m+1)}}}),$$
\p be the $m$-th Artin representation.
The Milnor invariants of length $\leq m$
vanish for a pure braid $\sigma \in P_n$,
if and only if all the longitudes 
of $\sigma$ are trivial in 
$F / F^{(m)}$; or equivalently
${\cal A}_m(\sigma)=1$.
It follows from the results in [BN2] (see also [H-M])
that ${\cal A}_m(\sigma)=1$ if and only if
$\sigma \in P_n^{(m)}$ or equivalently
the finite type invariants of orders $\leq m-1$
of $\sigma$ are the same as these of the
trivial braid in $P_n$.

For {\it pure links} (ordinary closures of pure braids)
a stronger conlusion than
Theorem 3.1 follows immediately from the discussion above.
However, most of $n$-component and  $n$-bridge
links are not  pure links. In particular,
many of the Brunnian examples  constructed by
``Bing doubling" in [C2],  are shown not to be pure
links. Moreover, one can see that if a link $L=L(\sigma)$,
whose ${\bar \mu}$-invariants
of length $\leq m$ vanish,
is the closure of a pure plat $\sigma$,
then we don't  necessarily have
${\cal A}_m(\sigma)=1$.
Thus, Theorem 1 doesn't in general follow from
the results in [BN2]. 
\medskip
b) By Corollary 6.10 of [H-M] the Kontsevich integral
of a (zero framed) link with vanishing
${\bar \mu}$-invariants of length $\leq n$,
can be written as a linear combination
of {\it Chinese character diagrams} 
that contain no tree components in degrees
$\leq n$. In view of our Theorem 1, it becomes interesting
to investigate the form that the  Kontsevich
integral takes on  closures of pure plats.
\medskip
\medskip
\p {\largebold References}
\medskip
\smallskip

\item {[B1]}  J. S. Birman: {\sl Braids, links and mapping
class group}, Ann. of Math. Studies, vol. {82},
Princeton University Press, 1974.

\item {[B2]} $\underline {\phantom {aaaaaaaaa}}$: {\sl  New  points of  view
  in knot theory}, Bulletin of AMS {\bf 28}(1993), 253--287.

\item {[BN1]}  D. Bar-Natan:
{ \sl On the Vassiliev knot invariants}, Topology {\bf 34}(1995), 423--472.

\item {[BN2]} $\underline {\phantom {aaaaaaaaaa}}$:
{\sl Vassiliev homotopy string link invariants},
J. of Knot Theory and Ramifications {\bf 4}(1995), 13--32. 

\item {[C1]} T. D. Cochran: {\sl Derivatives of Links:
Milnor's Concordance Invariants and 
Massey's Products,} Memoirs AMS {\bf 84}(1990).

\item {[C2]} $\underline {\phantom {aaaaaaaaaa}}$:
{\sl Non-trivial links and plats with
trivial Gassner matrices}, Math. Proc.Camb. Phil. Soc. {\bf 119}(1996), 
43--53.

\item {[C3]} $\underline {\phantom {aaaaaaaaaa}}$:
{\sl Links with trivial Alexander's module
but non-vanishing Massey products}, Topology {\bf 29}(1990), 189--204.

\item {[C-O]} $\underline {\phantom {aaaaaaaaaa}}$ and  K. Orr:
{\sl Not all links are concordant to boundary links},
Annals of Math. {\bf 138}(1993), 519--554.

\item {[G]} M. N. Gusarov:
{\sl On $n$-equivalence of knots and invariants of finite degree},
in {\sl Topology of manifolds and varieties}, ed. by O. Viro, 
Adv. Soviet Math., vol. {18}, AMS, 1994, 173--192.

\item {[H-M]} N. Habegger and G. Masbaum: {\sl
The Kontsevich integral and Milnor's Invariants},
preprint 1997.

\item {[K-L]} E. Kalfagianni and X.-S. Lin:
{\sl Regular Seifert surfaces and Vassiliev knot invariants},
preprint 1997 (available at {\tt http://www.math.ucr.edu/$^\sim$xl}).

%\item {[L]} T. Q. T. Le:
%{\it An invariant  of integral homology 3-spheres
%which is universal for all finite type invariants}, 
%Solitons, geometry, and topology, Amer. Math. Soc. Transl. Ser. 2,
%179, 
%AMS, Providence, RI, 1997.

\item{[Li]} X.-S. Lin: {\sl Power series expansions and invariants of links},
in {\sl Geometric Topology},
ed. by W. Kazez, Studies in Advanced Mathematics, vol. 2, AMS/IP, 1996, 
184--202.

\item {[M1]} J. Milnor: {\sl Link groups}.
Annals of Math. {\bf 2} (1954), 145-154.

\item {[M2]} J. Milnor:
{\sl Isotopy of links}, in {\sl Algebraic Geometry and Topology},
Princeton University Press, 1957, 280-386.

\item {[Ro]} D. Rolfsen: {\sl Knots and Links},
MLS, vol. 7, Publish or Perish, 1976.

\item{[V]} V. Vassiliev: {\sl
Cohomology of knot spaces}, in {\sl Theory of singularities and 
its applications}, ed. by V. I. Arnold,
Adv. Sov. Math, vol. 1, AMS, 1990, 23-68.

\medskip
\medskip
\p{\Smaller {\small D}EPARTMENT OF {\small M}ATHEMATICS,
 {\small R}UTGERS
{\small U}NIVERSITY, {\small N}EW {\small B}RUNSWICK, 
{\small NJ}
{\small 08903}}
\ \ \ \ \ \ \ \ {\vsmall  ekal@math.rutgers.edu}
%\medskip

\p{\small {Address as of September 1st 1998:}}
{\Smaller {\small D}EPARTMENT OF {\small M}ATHEMATICS,
 {\small M}ICHIGAN {\small S}TATE
{\small U}NIVERSITY, {\small E}AST {\small L}ANSING, 
{\small MI}
{\small 48824}}
\smallskip
\smallskip
\p{\Smaller {\small D}EPARTMENT OF {\small M}ATHEMATICS,
 {\small U}NIVERSITY OF
{\small C}ALIFORNIA, {\small R}IVERSIDE, 
{\small CA}
{\small 92521}}
\ \ \ \ \ \ \ \ {\vsmall  xl@math.ucr.edu}

\end